\documentclass[11pt, a4paper]{article}

\usepackage{amsmath,amssymb,epsfig,color}

\newcommand{\dd}{{\mathrm{d}}}
\newcommand{\ee}{{\mathrm{e}}}

\newcommand{\OO}{{\mathrm{O}}}

\newcommand{\sign}{{\mathrm{sign} \,}}

\newcommand{\ub}{\bar{u}}

\newcommand{\ut}{\tilde{u}}

\newcommand{\sigmab}{\bar{\sigma}}

\newcommand{\sigmat}{\tilde{\sigma}}

\newcommand{\bv}{{\boldsymbol{v}}}

\newcommand{\dif}[2]{\frac{\dd #1}{\dd #2}}
\newcommand{\pd}[2]{\frac{\partial #1}{\partial #2}}

\newcommand{\Dif}[2]{{\dd #1}/{\dd #2}}
\newcommand{\Pd}[2]{{\partial #1}/{\partial #2}}

\newcommand{\IN}[1]{{\, \dd #1}}

\newcommand{\m}{{\, \mbox{m}}}
\newcommand{\s}{{\, \mbox{s}}}

\newcommand{\K}{{\, \mbox{K}}}

\title{Transient Thermal Behaviour in a Model of Linear Friction Welding}
\author{A. A. Lacey\footnote{Corresponding author}
\\
\& \\
C. Voong\\
 \\
School of Mathematical and Computer Sciences,\\
and the Maxwell Institute for Mathematical Sciences,\\
Heriot-Watt University,\\
Edinburgh, EH14 4AS, UK}
\date{}

\begin{document}

\maketitle
\begin{abstract}

We derive a non-local model for the evolution of temperature
in workpieces being joined by linear friction welding.
The non-locality arises through the velocity being fixed
by the temperature gradient at the weld. Short- and long-time
behaviours are considered.

\end{abstract}

{\bf Keywords } Friction welding. Thermal modelling. Lubrication theory. Non-local problem.



\section{ Introduction} \label{intro} \setcounter{equation}{0}


We consider how a model of
linear friction welding predicts the way temperature approaches equilibrium.
In linear friction welding, two components,
usually both metal or metal alloy, are rubbed together to heat
their surfaces of contact, and then forced together while maintaining
a linear oscillatory motion; similar procedures are carried out in
the related methods of rotary and stir friction welding. (See \cite{LMY}
- \cite{WJ}.) A thin softened layer forms between the
components (the layer eventually makes up the weld after the process
has ceased).

If the oscillations and squeezing are maintained
long enough, a steady state can be reached, with periodic behaviour
in the thin, softened, welding region, and, to leading order, constant
temperature elsewhere. The unpublished Study Group report \cite{report}
and paper \cite{LV} looked at models for the temperature and
mechanical deformation in such steady states. To be more precise,
in \cite{report} the velocity towards the weld and the temperature
were constant while the velocity parallel to the weld was periodic
(with, to leading order, constant magnitude). In \cite{LV}, the
situation was either as in \cite{report}, or the velocity and
temperature could vary periodically in the thin welding region.

A key part of the model in both \cite{report} and \cite{LV} is that
virtually all the deformation and heat dissipation is confined to
the central, relatively warm, soft layer, where
the material can be modelled as a non-linear liquid and lubrication
theory can be applied.

In the paper \cite{KLNV} a transient model for the thermo-mechanical
behaviour of the thin central layer, where the welding occurs, was
considered. In particular, for a ``soft material'', for which the
temperature dependency in the stress-strain relation can be taken
to be exponential in the thin layer, it was shown that the non-local
PDE for temperature which applied in the thin layer had a global
solution, and that this solution tended to a steady state. Numerical
solutions for the corresponding non-local PDEs for ``hard material'',
for which the temperature dependence was locally a power,
indicated the same good behaviour. This variation with respect to
time was very fast, over a short time scale associated with thermal
diffusion on the scale of the width of the soft layer. The present
paper looks at much slower changes, associated with thermal
diffusion on the scale of the components, referred to as
``workpieces'', being welded.

\

We start by outlining some key parts of the modelling of the lubricating
layer as done in \cite{report} and \cite{LV}, as this is needed to build
our model for the time-dependent thermal problem. These are then used
to derive a key relation between the velocity of the rigid workpieces,
away from the softened layer, and the matching temperature gradient
between the rigid and softened regions. Because of the appearance
of a function of an end value of temperature gradient in the convective
term in the heat equation which applies in the rigid region,
this PDE is non-local. The qualitative
behaviour of solutions of these non-local PDEs for cases of both
hard and soft materials are briefly discussed, and one set of numerical
solutions is presented.


\section{The Model of the Lubricating Layer for a Hard Material} \label{model} \setcounter{equation}{0}

\begin{figure}[h]
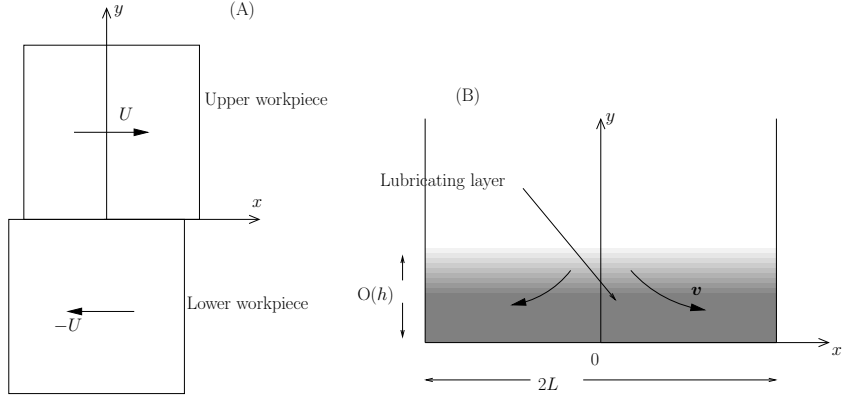

\begin{center}
\scalebox{0.33}{\input{1A.tex}}
\scalebox{0.33}{\input{1B.tex}}
\end{center}
\caption{(A) Upper workpiece sliding in $x$ direction with velocity
$U(t)$ and lower one with velocity $-U(t)$. (B) Upper workpiece
with its thin, lubricating, softened, welding layer, of width
O($h$), $h \ll L$, with $L$ the half-length of each workpiece.}\label{motion}
\end{figure}

The physical situation is indicated by Fig.~\ref{motion}, showing,
for simplicity, a two-dimensional symmetric case.

The workpieces have sides $x= \pm L$ which are taken to be effectively
insulated. Away from the thin softened layer, of width of order $h$,
where $h \ll L$, the workpieces move rigidly and linearly so that the
velocity $\bv = (u,v)$ is given by
\[
\bv \sim (U(t),0) \, \mbox{ for } \, y \gg h \, .
\]
Following \cite{report} and most of \cite{LV}, the simple case of
\[
U(t) = U_e \sign (\sin \omega t)
\]
is discussed in the present paper.

Focusing attention on the upper workpiece and times when $U>0$
(to avoid writing modulus signs), the basic model for the weld
in $0 < x < L$, $y > 0$ is here taken as:

\begin{equation} \label{eqncont}
\pd{u}{x} + \pd{v}{y} = 0 \, ,
\end{equation}
the usual continuity equation (the flow can be regarded as
incompressible);

\begin{equation}
\sigma = \kappa (T) \left( \pd uy \right)^{1/4} ,
\label{rheol}
\end{equation}
the constitutive law for shear stress $\sigma$ in terms of temperature
$T$ and strain rate $\partial u/\partial y$, with $u$ being the
component of velocity parallel to the weld (in the $x$ direction);

\begin{equation}\label{eqnmodtemp}
\kappa(T) = \kappa_m \left( 1 - \frac{T}{T_m} \right)
\exp \left( \frac{T_a}{T_m} { \left( \frac{T_m}{T} - 1 \right) } \right) ,
\end{equation}
the form of the temperature dependence in the constitutive law
--  here $\kappa_m$ is a material constant, $T_a$ is an
activation temperature and $T_m$ is the melting temperature;

\begin{equation}
\pd py = 0 \, , \quad \pd \sigma y = \pd px \, ,
\label{lub}
\end{equation}
the usual momentum-balance equations for two-dimensional lubrication
theory;

\begin{equation}\label{eqnheat}
k \pd{^2T}{y^2} + \sigma \pd ux = 0 \, ,
\end{equation}
the heat equation for $T$, neglecting any $x$ variation of temperature,
presently neglecting the time derivative of $T$ and the
convective term, and accounting for
viscous dissipation.

Note that in (\ref{eqnmodtemp}), the presence of the $(1 - T/T_m)$
factor models the stress becoming negligible as temperature approaches
melting point. The density $\rho$, specific heat $c_p$ and thermal
conductivity $k$ are all regarded as constant.

\

For a symmetric problem, the conditions
\begin{equation}
\pd px = 0 \, \mbox{ on } \, x = 0 \, , \quad
\pd vy = 0 \, \mbox{ on } \, y = 0 \, , \quad
\pd Ty = 0 \, \mbox{ on } \, y = 0 \, ,
\label{symmetry}
\end{equation}
are all taken to hold.

At the edge of the weld, pressure is atmospheric, so
\[
p = 0 \, \mbox{ at } \, x = L
\]
is imposed.

For distances large compared with the weld thickness, the velocity
component parallel to the weld is specified,
\[
u \to U_e \, \mbox{ as } \, y/h \to \infty \, ,
\]
while the normal component tends to some other value,
\[
v \to - V \, \mbox{ as } \, y/h \to \infty \, ,
\]
with the approach speed $V$ to be determined.

The workpieces are forced together with a specified average pressure
$P > 0$ so that
\[
\frac 1L \int_0^L p \IN x = P \, .
\]

The temperature should match with that in an outer region, where
the workpiece simply moves as a rigid body (and where there is no
dissipation term in the heat equation) so we shall require
\begin{equation}
\pd Ty \sim - G \, \mbox{ for } \, y/h \gg 1 \, ,
\label{heatmatch}
\end{equation}
with $G > 0$ to be found later.

\

Following \cite{report} (and \cite{LV}), the basic approach in solving the
problem is to decouple the velocity component $u$ and stress $\sigma$
into their ``sliding'' and ``squeezing'' parts,
\[
u \sim \ut + \ub \, \mbox{ and } \, \sigma \sim \sigmat + \sigmab \, ,
\quad \mbox{ where } \, \ub \ll \ut \, \mbox{ and } \, \sigmab \ll \sigmat
\, .
\]
The relative sizes of the sliding and squeezing parts can be justified,
\cite{LV}, in the limit of
\[
\epsilon^2 = \frac P{\rho c_p (T_m - T_e)} \to 0
\]
for $T_e$ the ambient temperature, say temperature at $y = l > 0$,
$l$ the length of the workpiece. Then, following \cite{LV}, asymptotic
approximations
\[
u \sim u_0 + \epsilon u_1 + \dots \, , \quad
\sigma \sim \sigma_0 + \epsilon \sigma_1 + \dots \, , \quad
T \sim T_0 + \epsilon T_1 + \dots
\]
can be sought. The leading terms $u_0$ and $\sigma_0$ can be identified
with $\ut$ and $\sigmat$ respectively, $\epsilon u_1$ with $\ub$ and
$\epsilon \sigma_1$ with $\sigmab$. Note that using these asymptotic
expansions gives somewhat different numerical values for quantities
such as approach speed $V$ from those found by the simpler
decoupling method. However, the orders of magnitude of these quantities,
and the way they depend on constants (such as $\kappa_m$
or $k$) or controllable physical parameters (such as $U_e$ or $P$)
are the same for the two approaches, \cite{LV}. Hence, for simplicity, only
the decoupling method is briefly described here.
The analysis covered in \cite{LV} is partially repeated as the outcomes
there are not given in the form required for later use in the present paper.

\

We first simplify the temperature dependency.
For a ``hard material'', to get appreciable softening, the temperature
in the thin welding layer must be very close to melting point
and the temperature
function (\ref{eqnmodtemp}) can then be approximated by
\begin{equation}
\kappa(T) \sim \kappa_m (1 - T/T_m) \, .
\label{templin}
\end{equation}

\

Turning to the decoupling approach,
the second momentum balance equation (\ref{lub}) separates as
\begin{equation}
\pd\sigmat y = 0 \, \mbox{ (sliding motion) } \qquad \mbox{ and }
\qquad \pd \sigmab y = \dif px\, \mbox{ (squeezing motion),}
\label{deco1}
\end{equation}
the constitutive law (\ref{rheol}), on using (\ref{templin}), as
\begin{equation}
\sigmat = \kappa_m (1 - T/T_m) \left( \pd\ut y \right)^{1/4} \qquad
\mbox{ and } \qquad \sigmab = \frac 14 \sigmat^{-3} \kappa_m^4
(1 - T/T_m)^4 \pd\ub y \, ,
\label{deco2}
\end{equation}
and the continuity equation (\ref{eqncont}) as
\begin{equation}
\pd\ut x = 0 \qquad \mbox{ and }
\qquad \pd \ub x + \pd vy = 0 \, ,
\label{deco3}
\end{equation}
while the thermal equation (\ref{eqnheat}) is approximately
\begin{equation}
k\dif{^2 T}{y^2} + \sigmat \pd\ut y = 0 \, .
\label{deco4}
\end{equation}

\

The matching conditions with the outer region are now
\begin{equation}
\ut \to U_e \, , \quad \ub \to 0 \, , \quad v \to - V  \quad
\mbox{ and } \quad \pd Ty \to - G \quad \mbox{ as } \quad \frac yh \to  \infty \, .
\label{deco5}
\end{equation}

The symmetry conditions to be imposed are
\begin{equation}
\ut = 0 \, , \quad v = 0 \, , \quad \pd\ub y = 0 \, , \,
\mbox{ and } \, \pd Ty = 0 \, \mbox{ on } \, y = 0
\label{deco6}
\end{equation}
and
\begin{equation}
\ub = 0 \, \mbox{ and } \, \dif px = 0 \, \mbox{ on } \, x = 0 \, .
\label{deco7}
\end{equation}

Finally,
\begin{equation}
p = 0 \, \mbox{ on } x = L \quad \mbox{ and } \quad
\frac 1L \int_0^L p \IN x = P \, .
\label{deco8}
\end{equation}

\

A solution for $\ub$ is sought in the form $\ub = xw(y)$, while,
because of the first of (\ref{deco1}), of (\ref{deco2}) and
of (\ref{deco3}), and because $T$ is independent of $x$, $\sigmat$
is constant so that $w(y)$ and $p(x)$ must satisfy
\[
\dif vy = - w \, , \qquad \dif px = \pd{}y \left(
\frac{\kappa_m^4}{4\sigma^3} \left( 1 - \frac T{T_m} \right)^4
\pd\ub y \right) =
\frac{\kappa_m^4}{4\sigma^3} \dif{}y \left(
\left( 1 - \frac T{T_m} \right)^4 \dif wy \right) x \, ,
\]
where, for convenience, $\sigmat$ has been replaced by $\sigma$
(they agree to leading order). Thus
\[
p = \frac{3P(L^2 - x^2)}{2L^2} \, ,
\]
and
\[
\frac{\kappa_m^4}{4\sigma^3} \dif{}y \left(
\left( 1 - \frac T{T_m} \right)^4 \dif wy \right)
+ \frac{3P}{L^2} = 0 \, ,
\]
with
\[
w \to 0 \mbox{ as }
\frac yh \to \infty \, \mbox{ and } \, \dif wy = 0 \mbox{ at } y = 0
\]
from the second of (\ref{deco5}) and the third of (\ref{deco6}). Then
\begin{equation}
w = \frac{12 P \sigma^3}{L^2 \kappa_m^4} \int_y^\infty
\frac{y \IN y}{(1 - T/T_m)^4} \, \mbox{ and } \,
V = \frac{12 P \sigma^3}{L^2 \kappa_m^4} \int_0^\infty
\frac{y^2 \IN y}{(1 - T/T_m)^4} \, .
\label{Veqn}
\end{equation}
Meanwhile, from the first of (\ref{deco2}), of (\ref{deco5})
and of (\ref{deco6}),
\begin{equation}
U_e = \frac{\sigma^4}{\kappa_m^4} \int_0^\infty
\frac{\IN y}{(1 - T/T_m)^4} \, ,
\label{deco9}
\end{equation}
while the heat equation (\ref{deco4}) becomes
\begin{equation}
k \dif{^2 T}{y^2} + \frac{\sigma^5}{\kappa_m^4 (1 - T/T_m)^4} = 0 \, ,
\label{heathard}
\end{equation}
which, with (\ref{deco9}) and the last of (\ref{deco5}) and of (\ref{deco6}),
yields
\[
G = \frac{\sigma^5}{k \kappa_m^4} \int_0^\infty
\frac{\IN y}{(1 - T/T_m)^4}
\]
\begin{equation}
\mbox{so } \qquad G = \frac{\sigma U_e}k \, .
\label{Geqn}
\end{equation}

\

Note that all these integrations to infinity, in (\ref{Veqn}) to
(\ref{Geqn}), mean taking $y/h \to \infty$ with the inner approximations
valid (in the soft layer). It is appropriate to scale distance
with $h$, $y = h\eta$, and also carry out other scalings,
$\sigma = Ss$, $T = T_m - \delta\theta$,
$V = V_EV^*$, $G = G_E G^*$,
where $h =$ typical width of lubricating layer, $S =$
size of stress, $\delta =$ size of temperature variation in layer,
$V_E =$ size of approach velocity and the scaling temperature
gradient $G_E$ is that taken at equilibrium
in a ``large'' workpiece. The non-dimensional problem (2.45) - (2.46)
of \cite{LV} almost results:
\begin{equation}
V^* = 12 s^3 \int_0^\infty \eta^2 \theta^{-4} \IN\eta \, \mbox{ and}
\label{Veqn1}
\end{equation}
\begin{equation}
\dif{^2 \theta}{\eta^2} = s^5 \theta^{-4} \, \mbox{ for } 0 < \eta < \infty \,
, \quad \dif\theta\eta = 0 \, \mbox { at } \, \eta = 0 \, , \quad
\dif\theta\eta \to G^* \, \mbox { as } \, \eta \to \infty \, ,
\label{thetaeqn}
\end{equation}
where $T_e$ is the ambient temperature of the workpiece, taken at
sufficiently large distances from the weld.
(The second part of (2.45) of \cite{LV}, $s = V^*$, is not required
here, as that was derived from matching with a steady temperature in
the far field.) The scaling constants satisfy
\[
\frac{k\delta}{h^2} = \frac{T_m^4 S^5}{\kappa_m^4 \delta^4} \, , \quad
\frac{\delta}{h} = G_E \, , \quad
G_E = \frac{S U_e}{k} \, ,
\]
\[
V_E = \frac{P T_m^4 S^3 h^3}{L^2 \kappa_m^4 \delta^4} \quad
\mbox{ and } \quad G_E = \frac{V_E \rho c_p (T_m - T_e)}{k} \, .
\]
(The first, third and fourth of these are got by having balances
in (\ref{heathard}), (\ref{Geqn}) and (\ref{Veqn}) respectively,
while the second is obtained from the matching condition for
temperature gradient, (\ref{deco5}).
The last comes from having $T$ satisfy a steady convection-diffusion
equation, with velocity $-V_E$, in the outer, rigid region, with $T
\to T_e$ as $y \to \infty$.)
After some manipulation,
\[
V_E = (kT_m/\kappa_m)^{4/3} (P^{1/2}/L) (\rho c_p (T_m - T_e))^{-1/2}
U_e^{-2/3} \, ,
\]
\[
S = (kT_m/\kappa_m)^{4/3} (P^{1/2}/L) (\rho c_p (T_m - T_e))^{1/2}
U_e^{-5/3} \, ,
\]
\[
h = (k T_m / \kappa_m)^{4/3} U_e^{-5/3},
\]
\begin{equation}
\delta = k^{5/3}(T_m / \kappa_m)^{8/3} (\rho c_p (T_m-T_e))^{1/2}
(P^{1/2}L^{-1}) \, U_e^{-7/3},
\label{eqndelta}
\end{equation}
\[ \mbox{ and } \,
G_E = k^{1/3} (T_m/\kappa_m)^{4/3} (P^{1/2}/L)
(\rho c_p (T_m - T_e))^{1/2} U_e^{-2/3} \, .
\]

Eqn.~(\ref{Geqn}) scales to
\[ G^* = s \]
and then the change of variable $\theta = s \varphi$ alters
(\ref{thetaeqn}) to
\begin{equation}
\dif{^2 \varphi}{\eta^2} = \varphi^{-4} \, \mbox{ for } 0 < \eta < \infty \,
, \quad \dif\varphi\eta = 0 \, \mbox { at } \, \eta = 0 \, , \quad
\dif\varphi\eta \to 1 \, \mbox { as } \, \eta \to \infty
\label{thetan}
\end{equation}
and (\ref{Veqn1}) to
\begin{equation}
G^* V^* = sV^* = 12 \int_0^\infty \eta^2 \varphi^{-4} \IN\eta = N \, ,
\label{Veqn2}
\end{equation}
where $N$ is a numerical constant. From this decoupling approach of
\cite{report} and \cite{LV}, $N \approx 8.123$. 
Returning to dimensional quantities,
\begin{equation}
VG = N k^{5/3} (T_m/\kappa_m)^{8/3} (P/L^2) U_e^{-4/3} \, .
\label{VG}
\end{equation}
Note that the factor $\rho c_p (T_m - T_e)$, which comes from
considering long-time behaviour and might suggest influence at a
distance, has dropped out.

Note also that a different value of $N$ would be obtained from
using a formal asymptotic approach, as in Section~4 of \cite{LV}.


\section{The Non-Local Model for a Hard Material} \label{nonlocal} \setcounter{equation}{0}


In the outer region,
away from the welding layer, the workpiece is simply modelled with
regard to its thermal behaviour, regarding it as a rigid body moving
with speed $V(t)$ towards the weld (the oscillatory motion in the
$x$ direction can be disregarded):
\begin{equation}
\rho c_p \left( \pd Tt - V \pd Ty \right) = k \pd{^2 T}{y^2} \quad
\mbox{ for } \quad  0 < y < l \, , \, t > 0 \, .
\label{TPDE}
\end{equation}
We assume the side $x = L$ to be thermally insulated
so that $T = T(y,t)$. The temperature is taken to be ambient at the
far end of the workpiece and initially:
\begin{equation}
T = T_e \quad \mbox {for} \quad y = l \, , \, t \ge 0 \, ; \quad
T = T_e \quad \mbox {for} \quad 0 < y \le l \, , \, t = 0 \, .
\label{Tconds}
\end{equation}
As we are considering the case of a hard material, the weld should
be at melting temperature, so, to leading order,
\begin{equation}
T = T_m \quad \mbox {for} \quad y = 0 \, , \, t > 0 \, .
\label{Tcond}
\end{equation}
The problem is completed by fixing $V(t)$ through the non-local
condition (\ref{VG}):
\begin{equation}
V(t) = - M \left/ \pd Ty \right|_{y=0} \, ,
\label{Vnonloc}
\end{equation}
with
\begin{equation}
M = N k^{5/3} (T_m/\kappa_m)^{8/3} (P/L^2) U_e^{-4/3} \, .
\label{Mdefn}
\end{equation}
It is possible to non-dimensionalise the problem (as for the case of
soft materials, Sec.~\ref{soft}) but as no great simplification of
the problem results in this case, this is not done here.

\paragraph{Long-time solution.} We expect $T$ to tend towards the
steady-state, which satisfies
\[
k\dif{^2 T}{y^2} + \rho c_p V_\infty \dif Ty = 0 \, \mbox{ for } \,
0 < y < l \, , \quad T = T_m \, \mbox{ at } \, y = 0 \, , \quad
T = T_e \, \mbox{ at } \, y = l \, ,
\]
\[
\mbox{ with }
V_\infty = - M \Big/\dif Ty \, \mbox{ at } \, y = 0 \, .
\]
With large enough $l$ (so that $\exp(- \rho c_p V_\infty l/k) \ll 1$), \\
$T \sim T_e + (T_m - T_e)\exp(- \rho c_p V_\infty y/k)$, 
$\Dif Ty \sim - \rho c_p V_\infty (T_m - T_e)/k$ at $y = 0$, and
\[
V_\infty \sim \sqrt{ \frac{Mk}{\rho c_p(T_m - T_e)} }
\]
is the steady velocity of approach.

\paragraph{Short-time solution.} For a sufficiently short time, the
``hard-material approximation'', $\kappa(T) \sim \kappa_m(1 - T/T_m)$,
would not be valid and in particularly extreme cases we would have
$T \sim T_c(t)$, with $T_c$ not close to $T_m$, in the softened
layer. Such cases will be included in the discussions of short-time
regimes for soft materials in Sec.~\ref{soft}. For the present, we
take time to be small compared with the time scale of (\ref{TPDE})
({\it i.e.} $t \ll l^2 \rho c_p/k$, or $t \ll l_\infty^2 \rho c_p/k
= (T_m - T_e)/M$ if the length scale for the steady state,
$l_\infty = k/(\rho c_p V_\infty) = \sqrt{k(T_m - T_e)/(\rho c_p M)}$,
is small compared with the length of $l$ of the workpiece), but still
with $T \sim T_m$ near $y = 0$ and $\kappa(T) \sim \kappa_m(1 - T/T_m)$
in the welding layer, so that (\ref{TPDE}) - (\ref{Vnonloc}) apply.

For short times, high temperature gradients are expected near $y = 0$
so $V$ should be small and (\ref{TPDE}) reduces to the heat equation:
\begin{equation}
\rho c_p \pd Tt = k \pd{^2 T}{y^2} \quad 
\mbox{ for } \quad  0 < y < l \, , \, t > 0 \, .
\label{Theat}
\end{equation}
Combined with the boundary condition (\ref{Tcond}) and initial condition
$T = T_e$, this leads to the standard local similarity solution
\[
T \sim T_e + (T_m - T_e)\mathrm{erfc} \left( \frac x2 \sqrt{
\frac{\rho c_p}{k t}} \right) = T_e + \frac 2{\sqrt{\pi}} (T_m - T_e)
\int_{\frac x2 \sqrt{ \frac{\rho c_p}{k t}}}^\infty \ee^{-z^2} \IN z
\]
(see, for example, \cite{CJ} or \cite{OHLM}). In particular,
\[
- \left. \dif Ty \right|_{y=0} \sim (T_m - T_e) \sqrt{ \frac {\rho c_p}
{\pi k t}}
\quad \mbox{ so } \quad V \sim
\frac{M}{(T_m - T_e)} \sqrt{ \frac {\pi k t}{\rho c_p}} \, .
\]

\begin{figure}[h]
\begin{center}
\psfig{figure=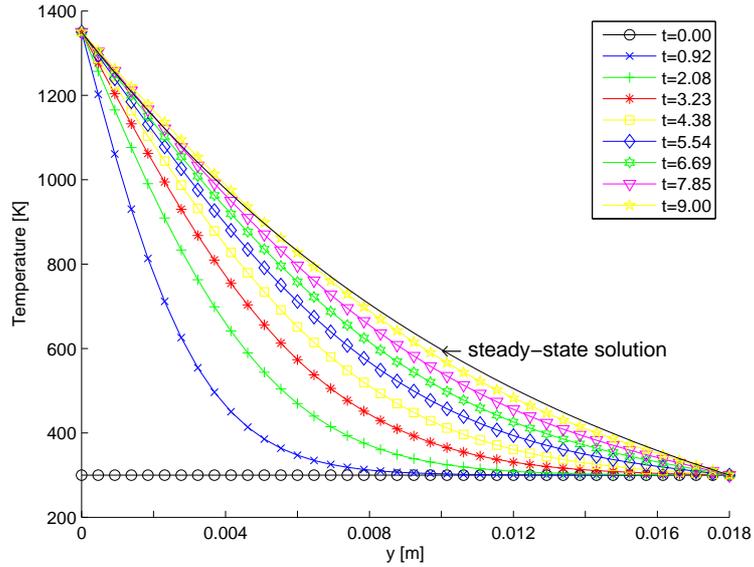,width=11cm}
\end{center}
\caption{Plot of temperature, $T$, for different times, taking
the case of a hard material, so that approximation (\ref{templin})
applies.}\label{numsol}
\end{figure}

\paragraph{Numerical solution.} The set of simulations carried out here
used a backward-time, central-space, implicit approximation.
Fig.~\ref{numsol} shows results for $k/(\rho c_p) = 4.52 \times 10^{-6}
\m^2 \s^{-1}$, $M = 28.024 \K \s^{-1}$, $l = 0.018 \m$, $T_e = 300 \K$,
$T_m = 1350 \K$ (so $l_\infty \approx 0.013\m \approx l$). Note the
high initial slope near $y = 0$ and approach, with increasing $t$,
towards the steady state.


\section{Soft Materials} \label{soft} \setcounter{equation}{0}


If the material is ``soft'' (or is hard but the temperature is well
below melting in the thin lubricating layer), the temperature
dependency function in the thin layer is now approximated by
\begin{equation}
\kappa(T) \sim K_c \exp \left( \frac{T_a^2}{T_c}(T_c - T) \right)
\label{softap1}
\end{equation}
where
\begin{equation}
K_c = \kappa_m \left( 1 - \frac{T_c}{T_m} \right)
\exp \left( \frac{T_a}{T_c} - \frac{T_a}{T_m} \right) \, ,
\label{softap2}
\end{equation}
with $T_c = T|_{y=0}$ generally varying with time.

Following the calculations of \cite{LV}, which means, for the simple
decoupling approach, replacing the negative fourth power of temperature
in the heat equation and integrals by an exponential,
so that (\ref{heathard}) is replaced by
\begin{equation}
k\dif{^2 T}{y^2} + \frac{\sigma^5}{K_c^4} \ee^{-4T_a(T_c - T)/T_c^2} = 0
\label{heatsoft}
\end{equation}
(see (3.6) in \cite{LV} where the dimensionless temperature is now
$\theta = T_a(T_c - T)/T_c^2$) and (\ref{Veqn}) by
\begin{equation}
V = \frac{12 P\sigma^3}{L^2 K_c^4} \int_0^\infty y^2 \exp \left(
\frac{-4T_a}{T_c^2} (T_c - T) \right) \IN y \, .
\label{Veq}
\end{equation}
The solution of (\ref{heatsoft}) subject to (\ref{heatmatch}) and
$\Pd Ty = 0$ on $y = 0$, where $T = T_c$, is
\begin{equation}
T = T_c - \frac{T_c^2}{2T_a} \ln \left( \cosh by \right) \, ,
\label{Tsol}
\end{equation}
\begin{equation}
\mbox{with } \, b = 2T_a G/T_c^2 \, \mbox{ and } \, b^2 =
2\sigma^5 T_a/(k K_c^4 T_c^2) \, .
\label{bs}
\end{equation}
Using (\ref{Tsol}) and the first of (\ref{bs}) in (\ref{Veq}) gives
\[
V = \frac{\pi^2 T_c^6 P \sigma^3}{8 K_c^4 T_a^3 L^2 G^3}
\]
(see \cite{LV}) while (\ref{bs}) combined with (\ref{Geqn}) (which holds
for any approximation of $\kappa(T)$) leads to
\begin{equation}
T_c^2 G^3/K_c^4 = 2T_a U_e^5/k^4 \, ,
\label{10D}
\end{equation}
and hence (\ref{Geqn}) gives
\begin{equation}
V = \frac{\pi^2 k^3 T_c^6 P}{8 K_c^4 T_a^3 L^2 U_e^3} =
\frac{\pi^2 P U_e^2 T_c^4}{4 L^2 k T_a^2 G^3} \, .
\label{10B}
\end{equation}
As $K_c$ is given in terms of $T_c$ by (\ref{softap2}), (\ref{10D})
can be written as
\begin{equation}
\kappa_m^{-4} \left( 1 - \frac{T_c}{T_m} \right)^{-4} G^3 T_c^2
\exp \left( \frac{4T_a}{T_m} - \frac{4T_a}{T_c} \right) =
\frac{2 T_a U_e^5}{k^4} \, .
\label{10E}
\end{equation}
Eqn.~(\ref{10E}), applied as a matching condition for the outer problem,
gives a boundary condition at $y=0$, while (\ref{10B}) now gives a
non-local effect through having a velocity for the outer problem
of the form
\begin{equation}
V = \frac{N T_c^4 P U_e^2}{k T_a L^2 G^3} \quad \mbox{ with }
\quad G = - \left. \pd Ty \right|_{y=0} \, .
\label{11Z}
\end{equation}
The value of $N$ according to the simple-minded decoupling approach
of \cite{report} and \cite{LV} outlined in Sec.~\ref{model} is
$\pi^2/4$, but a more complete asymptotic solution (see \cite{LV}) leads
to $N = \frac 34 ( \frac 52 - \frac {\pi^2}{12} )$.

\

The outer problem, representing just heat conduction and convection, is
again (\ref{TPDE}) with (\ref{Tconds}),
\begin{equation}
\rho c_p \left( \pd Tt - V \pd Ty \right) = k \pd{^2 T}{y^2} \quad 
\mbox{ for } \quad  0 < y < l \, , \, t > 0 \, ,
\label{TPDEa}
\end{equation}
\begin{equation}
T = T_e \quad \mbox {for} \quad y = l \, , \, t \ge 0 \, ; \quad
T = T_e \quad \mbox {for} \quad 0 < y \le l \, , \, t = 0 \, ,
\label{Tcondsa}
\end{equation}
but now with the second boundary condition on $y = 0$
based on (\ref{10E}),
\begin{equation}
\frac{T^2}{\kappa_m^{4} (1 - T/{T_m})^{4}} \left( - \pd Ty \right)^3
\exp \left( \frac{4T_a}{T_m} - \frac{4T_a}{T} \right) =
\frac{2 T_a U_e^5}{k^4}
\label{Tcondsc}
\end{equation}
and convective velocity given by
\begin{equation}
V = M(T) \left/ \left( - \pd Ty \right|_{y=0} \right)^3 \, ,
\label{Vnonloca}
\end{equation}
\begin{equation}
\mbox {for} \quad  M(T) = \frac{N P U_e^2}{k T_a L^2} \left( T|_{y=0}
\right)^4 \, .
\label{Tcondsb}
\end{equation}

\

The steady problem is
\[
k \dif{^2 T}{y^2} + \rho c_p V \dif Ty = 0 \, \mbox{ for } \, 0 < y < l \, ,
\quad \mbox{ with } T = T_e \, \mbox{ on } \, y = l \, ,
\]
\begin{equation}
\frac{T^2}{\kappa_m^{4} (1 - T/{T_m})^{4}} \left( - \dif Ty \right)^3
\exp \left( \frac{4T_a}{T_m} - \frac{4T_a}{T} \right) =
\frac{2 T_a U_e^5}{k^4} \, \mbox{ on } \, y = 0 \, ,
\label{Tsteadyb1}
\end{equation}
and
\begin{equation}
V = M \left/ \left( - \pd Ty \right|_{y=0} \right)^3 \, .
\label{Tsteadyb2}
\end{equation}
This fixes a limiting central temperature $T_\infty = T(0)$,
velocity $V_\infty$, length scale $l_\infty = k/(\rho c_p V)$ (which is
significant if it is less than $l$) and hence a corresponding time scale
for the problem (\ref{TPDEa}) - (\ref{Vnonloca}), $t_\infty =
l_\infty/V_\infty = l_\infty^2 \rho c_p/k$.

The dimensionless version of (\ref{TPDEa}) - (\ref{Vnonloca}),
given by scaling:
\[
y = l_\infty z \, , \quad t = t_\infty s \, , \quad
T = T_e + (T_\infty - T_e)\phi \, , \quad V = V_\infty V^* \, ;
\]
is then, assuming, for simplicity, that $l^* = l/l_\infty \gg 1$:
\[
\pd\phi s - V^*\pd\phi z = \pd{^2 \phi}{z^2} \quad 
\mbox{ for } \quad  z > 0 \, , \, s > 0 \, ;
\]
\[
\phi = 0 \quad \mbox {for} \quad s = 0 \, , \, z \ge 0 \, ;
\]
\begin{equation} \begin{array}{c} \displaystyle
\left( \frac{T_e + (T_\infty - T_e)\phi}{T_\infty} \right)^2
\left( \frac{1 - T_\infty/T_m}{1 - (T_e + (T_\infty - T_e)\phi)/T_m} \right)^4
\left( - \pd\phi z \right)^3 \\
\displaystyle = \exp \left( \frac{4T_a}{T_e + (T_\infty - T_e)\phi} -
\frac{4T_a}{T_\infty} \right)
 \, \mbox{ at } \, z = 0 \, ; \end{array}
\label{phi}
\end{equation}
\[
\mbox{and } \qquad V^* =
\left( \frac{T_e + (T_\infty - T_e)\phi|_{y=0}}{T_\infty} \right)^4
\left/ \left( - \pd\phi z \right|_{z=0} \right)^3 \, .
\]
Taking $T_\infty$ not close to $T_e$ with $T_\infty = \OO (T_e )$,
activation temperature to be large, $T_a \gg T_\infty$,
and assuming that $|\Pd\phi z|$ is not (exponentially) large, (\ref{phi}) reduces to
\[
\phi = 1 \, \mbox{ at } \, z = 0 \, ;
\]
the (dimensional) temperature near the centre of the weld varies by
only $\OO(T_\infty^2/T_a) \ll T_\infty$. With the assumptions
on how $T_e$ and $T_\infty$ compare, this means that the central
temperature varies much less than $T_\infty -T_e$.
The same will hold true for $l^*$ not large.

\paragraph*{Long-time solution.} Again we simply expect an approach to the
equilibrium, as determined above.

\paragraph*{Short-time solution.} Consideration of the model (\ref{TPDEa})
- (\ref{Vnonloca}) indicates a number of different regimes. In practice,
because of the high activation temperature $T_a$, these would last
exponentially short times, quite possibly comparable with, or smaller than,
the period of oscillation. In such extreme situations the model would
not be valid.

\paragraph*{(i) Initial stage: $T_c(t) = T|_{y=0}$ very close to $T_e$.}

From (\ref{Tcondsc}),
\[ \begin{array}{rl}
G \sim G_e = & \displaystyle \left( 2 T_a \kappa_m^4 k^{-4}
U_e^5 (1 - T_e/T_m)^4 T_e^{-2}
\exp \left( \frac{4T_a}{T_e} - \frac{4T_a}{T_m} \right) \right)^{1/3} \\
= &  \displaystyle G_0 \exp \left( \frac 43
\left( \frac{T_a}{T_e} - \frac{T_a}{T_m}
\right) \right) \, , \end{array}
\]
for
\begin{equation}
G_0 = \left( 2 T_a \kappa_m^4 k^{-4} U_e^5 (1 - T_e/T_m)^4 T_e^{-2}
\right)^{1/3} ,
\label{G0def}
\end{equation}
so $G$ is expected to be exponentially large for cases of interest, $V$ is
neglected and $T$ is given by
\[
\rho c_p \pd Tt \sim k \pd{^2 T}{y^2} \, , \, \mbox{ with } \,
- \pd Ty \sim G_e \mbox{ at } y = 0 \, .
\]
Hence, near $y = 0$, $T$ is given by a similarity solution,
\[
T \sim T_e + a t^{1/2} f(\eta) \, \mbox{ with } \, \eta = \left(
\frac{\rho c_p}{kt}\right)^{1/2} y \, , \quad a =
\left( \frac k{\rho c_p} \right)^{1/2} G_e
\]
and where $f$ satisfies
\[
\dif{^2 f}{\eta^2} + \frac\eta 2 \dif f\eta - \frac 12 f = 0 \,
\mbox{ for } \, \eta > 0 \, , \quad \mbox{ with } \, \dif f\eta = - 1
\, \mbox{ at } \, \eta = 0 \, .
\]
(For completeness, $\displaystyle f(\eta) = \frac 2{\sqrt{\pi}}
\left( \ee^{-\eta^2/4} - \frac\eta 2 \int_\eta^\infty \ee^{-\xi^2/4}
\IN\xi\right)$.)

This approximation for $T$ remains valid only as long as $G$ can be
approximated by $G_e$, which requires that $|T - T_e| \ll T_e^2/T_a$.
Hence this initial stage applies for
\[
t \ll \frac{\rho c_p}k G_0^{-2} \frac{T_e^4}{T_a^2} \exp \left(
\frac 83 \left( \frac{T_a}{T_m} - \frac{T_a}{T_e} \right) \right) \, .
\]

\paragraph*{(ii) First intermediate stage: $T_c(t) - T_e$
of order $T_e^2/T_a$.}

On writing $T = T_e + (T_e^2/T_a)\theta$,
\[
G \sim G_0 \exp \left( \frac 43 \left( \frac{T_a}{T_e} - \frac{T_a}{T_m}
\right) \right) \ee^{- \frac 43 \theta} = G_e \ee^{- \frac 43 \theta}
\]
so that $V^*$ (which is of size $G^{-3}$) is still exponentially
small and therefore negligible,
\[
\rho c_p \pd\theta t \sim k \pd{^2 \theta}{y^2} \, , \, \mbox{ with } \,
- \pd\theta y \sim \frac{G_e T_a}{T_e^2} \ee^{- \frac 43 \theta}
\mbox{ at } y = 0
\]
and $\theta = 0$ at $t = 0$.

Rescaling time, distance and temperature,
\begin{equation} \begin{array}{c} \displaystyle
t = \frac 9{16} \frac{\rho c_p}k \frac{T_e^4}{T_a^2} \frac 1{G_0^2}
\exp \left( \frac {16}9 \left( \frac{T_a}{T_m} - \frac{T_a}{T_e}
\right) \right) \tau \, , \\
\displaystyle y = \frac 34 \frac{T_e^2}{T_a} \frac 1{G_0}
\exp \left( \frac 43 \left( \frac{T_a}{T_m} - \frac{T_a}{T_e}
\right) \right) Y \, , \quad \theta = \frac 34 \varphi \, , \end{array}
\label{shortscale}
\end{equation}
leads to
\begin{equation}
\pd\varphi \tau \sim \pd{^2 \varphi}{Y^2} \, , \, \mbox{ with } \,
- \pd\varphi Y \sim \ee^{- \varphi}
\mbox{ at } Y = 0 \, ,
\label{shortPDE}
\end{equation}
and $\varphi = 0$ at $\tau = 0$. The form of the boundary condition
in (\ref{shortPDE}) suggests that there is no (exact) similarity
solution in this regime. (The decaying exponential in the boundary
condition distinguishes this phase of the problem from the very
early stage, when a constant condition applied  --  the constant
value being needed for the similarity solution.)

This regime will be left when $\theta$ and $\varphi$ are large which,
from (\ref{shortPDE}), clearly corresponds to $\tau \gg 1$, {\it i.e.}
\[
t \gg \frac{\rho c_p}k \frac{T_e^4}{T_a^2} \frac 1{G_0^2}
\exp \left( \frac {16}9 \left( \frac{T_a}{T_m} - \frac{T_a}{T_e}
\right) \right) \, .
\]

\paragraph*{(iii) Second intermediate stage: $T_c(t) - T_e \gg T_e^2/T_a$
but $T_\infty - T_c(t) \gg T_\infty^2/T_a$.} ($T_\infty$ being the
equilibrium temperature at the centre of the weld.)

The second condition states that we do not yet have $T_c$ close to
$T_\infty$. The velocity $V$ is still small so that the PDE can be
approximated by the heat equation but no further simplifications are
apparent; the condition at $y = 0$ is
\[
- \pd Ty \sim \left ( \frac{2 \kappa_m^4 T_a U_e^5 (1 - T/T_m)}
{k^{4} T^{2}} \right)^{1/3}
\exp \left( \frac 43 \left( \frac{T_a}{T_e} - \frac{T_a}{T_m} \right) \right)
\, .
\]

Note that each of these early stages could also be applicable to a
material regarded as hard  --  hard in that the steady maximum temperature has
to be so close to the melting temperature that the temperature-dependency
function $\kappa(T)$ is, locally, approximately linear. The second
intermediate stage would cease for $T_m - T_c$ of order $T_m^2/T_a$
in such cases. (There would then be time regimes, before $T_m - T_c$
falls to $\OO(\delta) = \OO(k^{5/3} (T_m/\kappa_m)^{8/3} (\rho c_p
(T_m - T_e))^{1/2} (P^{1/2} L^{-1}) U_e^{-7/3})$ (see (\ref{eqndelta}))
where $T_m - T_c$ and the matching $-\Pd Ty$
are so large that $\kappa(T)$ cannot be approximated by the linear
function, and $V(t)$ is not given by (\ref{Vnonloc}). However, during
such later intermediate times, $V(t)$ can still be regarded as
negligible and, for the
outer problem, $T = T_m$, to leading order, at $y = 0$.

\paragraph*{Numerical solution.} Avoiding the short-time regimes,
with associated exponentially short time scales (for large $T_a$),
the numerical method used for hard materials, Sec.~\ref{nonlocal},
can be employed again. Unsurprisingly, the results are qualitatively
very similar and are therefore not presented here.


\section{Conclusions} \label{conclude} \setcounter{equation}{0}


The steady-temperature modelling of \cite{report} and \cite{LV} has
been extended to consider the transient conditions which can occur
when linear friction welding is only run for finite times. In both
limiting cases of soft and hard materials, in which the
stress-strain-temperature constitutive law can be approximated in
different ways, the key problem to be solved here is an outer problem
for temperature which contains a non-local term. The non-locality
arises in the convective term in the heat equation, through the
velocity being fixed (thanks to the behaviour of the inner problem)
by the temperature gradient in the vicinity of the weld.

It can be observed that for an early enough stage in the process, so
that what might be regarded as a ``hard material'' (as its equilibrium
temperature gets so close to melting) has everywhere relatively
low temperatures, the approximate theory is the
same as that used for soft materials.
However, whether the workpieces are hard or soft, these early stages of
the process, with the central temperature significantly lower than its
equilibrium level, happen at rates exponentially quickly compared with
the later stages (see, for example, (\ref{shortscale})), assuming
that the activation temperature is large, $T_a \gg T_c(t)$. (For
$T_a = \OO(T_c(t))$, the exponential approximation (\ref{softap1})
used in analysing the thin welding layer becomes invalid.) For $T_c$
too small, the predicted changes will occur over time scales
comparable with, or even shorter than, the period of oscillation of
the workpieces, and again modelling assumptions needed for the analysis
of the thin layer would be violated. It should be emphasised that,
in practice, a preliminary stage of linear friction welding entails the
workpieces being rubbed against each other without them being forced
together. In this stage, simple surface friction results in localised
heating at the contact interfaces, without material deformation. This
present paper only considers what happens subsequently: when, after
such initial preheating the workpieces are pressed together, so that
material deformation starts, while continuing with the oscillatory
motion.



\end{document}